\documentclass[a4paper,10pt]{amsart}
\usepackage{amssymb}

\numberwithin{equation}{section}
\newtheorem{theorem}{Theorem}[section]
\newtheorem{prop}[theorem]{Proposition}
\newtheorem{lem}[theorem]{Lemma}

\theoremstyle{remark}

\newcommand{\R}{\mathbb{R}}

\author[J.~Benameur]{Jamel Benameur}
\author[H.~Orf]{Hajer Orf}
\address{I.S.S.A.T. Gab\`es, University of Gab\`es\\
Gab\`es, Tunisia}
\email{\sl jamelbenameur@gmail.com\;and\;jlbenameur@hotmail.com}
\address{Department of Mathematics, Faculty of Science, Gab\`es\\
Gab\`es, Tunisia}
\email{\sl hajerorf17@gmail.com}

\title[On the Blow up criterion of 3D-Navier-Stokes equation in $\dot H^{5/2}$]
{On the Blow up criterion of 3D-Navier-Stokes equation in $\dot H^{5/2}$}

\begin{document}
\begin{abstract}
In this paper, we prove two results about the blow up criterion of the three-dimensional incompressible Navier-Stokes equation in the sobolev space $\dot H^{5/2}$. The first one improves the result of \cite{CZ}. The second deals with the relationship of the blow up in $\dot H^{5/2}$ and some critical spaces. Fourier analysis and standard techniques are used.
\end{abstract}


\maketitle
\tableofcontents


\section{Introduction}

The $3D$ incompressible Navier-Stokes equations are given by:
$$
\left\{
  \begin{array}{lll}
     \partial_t u
 -\nu\Delta u+ u.\nabla u  & =&\;\;-\nabla p\hbox{ in } \mathbb R^+\times \mathbb R^3\\
     {\rm div}\, u &=& 0 \hbox{ in } \mathbb R^+\times \mathbb R^3\\
    u(0,x) &=&u^0(x) \;\;\hbox{ in }\mathbb R^3,
  \end{array}
\right.
\leqno{(NSE)}
$$
where $\nu>0$ is the viscosity of fluid, $u=u(t,x)=(u_1,u_2,u_3)$ and $p=p(t,x)$ denote respectively the unknown velocity and the unknown pressure of the fluid at the point $(t,x)\in \mathbb R^+\times \mathbb R^3$, and $(u.\nabla u):=u_1\partial_1 u+u_2\partial_2 u+u_3\partial_3u$, while $u^0=(u_1^o(x),u_2^o(x),u_3^o(x))$ is an initial given velocity. If $u^0$ is quite regular, the divergence free condition determines the pressure $p$.\\

In 1934, Leray \cite{JL02} showed that there exists an  absolute constant c such that, if  $\|u(t)\|_{\dot{H}^{5/2}}$ continuous on $[0,T^*[,$ and blow up at time $T^*,$ then
$$\|u(t)\|_{\dot{H}^1}\geq C(T^*-t)^{-1/4}.$$
Morever, he stated the bound for $L^p$ norms for $3<p<\infty$ (without proof) $$\|u(t)\|_{L^p}\geq C(T^*-t)^{\frac{3-p}{2p}}.$$
More recently, there have been a number of papers that treat the problem of blow up in Sobolev spaces $\dot{H}^s$ $s>1/2.$\\
Benameur (2010) \cite{JB1} showed that for $s>5/2$, $$\|u(t)\|_{\dot{H}^s}\geq C(T^*-t)^{\frac{-s}{3}},$$
which was improved by Robinson, Sadwiski, Silva (2012) \cite{4}
$$\|u(t)\|_{\dot{H}^s}\geq \begin{cases}
                                C(T^*-t)^{-(2s-1)/4}\quad 1/2<s<5/2,\; s\neq 3/2,\\
                                C(T^*-t)^{-s/5}\quad s>5/2.
                              \end{cases}$$
The border cases $s=3/2$ and $s=5/2$ required separate treatment. In \cite{6},
, Cortissoz, Montero, and Pinilla (2014)  proved lower bounds in $\dot H^{3/2}$ and $\dot H^{5/2}$ at the optimal rates but with logarithmic corrections,
$$\|u(t)\|_{\dot H^{3/2}}\geq\frac{c}{\sqrt{(T^* - t)|\log(T^* - t)|}},\;
\|u(t)\|_{\dot H^{5/2}}\geq\frac{c'}{(T^* - t)|\log(T^* - t)|},$$
where in both cases c depends on $\|u_0\|_{L^2}$. Recently,  in \cite{5} the authors  proved
$$\|u(t)\|_{\dot H^{3/2}}\geq\frac{c}{\sqrt{(T^* - t)}},$$
which we refer to as a "strong blowup estimate," and
$$\limsup_{t\rightarrow T^*}(T^* - t)\|u(t)\|_{\dot H^{5/2}}\geq c,$$
which we refer to as a "weak blowup estimate." They also show a strong blowup
estimate in the Besov space $\dot{B}^{5/2}_{2,1}$, which has the same scaling as $\dot{H}^{5/2},$
$$\|u(t)\|_{\dot{B}^{5/2}_{2,1}}\geq C(T^*-t)^{-1}.$$
The interesting and open question is the strong blow up estimate
$$\|u(t)\|_{\dot{H}^{5/2}}\geq C(T^*-t)^{-1}.$$
The some kind of question appears for Lei-Lin espace
$$\|u(t)\|_{\chi^1}\geq C(T^*-t)^{-1},$$
and, $$\|\nabla u(t)\|_{L^\infty}\geq C(T^*-t)^{-1}$$
i.e., a bound with the optimal rate in a space with the same scaling as $\dot{H}^{5/2}$.

Our first result is the following:
\begin{theorem}\label{theo1}
Let $u\in{\mathcal C}([0,T^*),H^{5/2}(\mathbb{R}^3)$ be a maximal solution of Navier-Stokes system. If, $T^*$ is finite then, there is a constant $c_0=c_0(\nu,\|u^0\|_{L^2})$ such that
\begin{eqnarray}\label{enq111}
\liminf_{t\nearrow T^*}(T^*-t)\sqrt{|\ln(T^*-t)|}\|u(t)\|_{\dot H^{5/2}}\geq c_0.
\end{eqnarray}
\end{theorem}
The second result is the following
\begin{theorem}\label{theo2}
Let $u\in{\mathcal C}([0,T^*),H^{5/2}(\mathbb{R}^3)$ be a maximal solution of Navier-Stokes system. If, $T^*$ is finite then, there is a universal constant $c_1>0$ such that
\begin{eqnarray}\label{enq112}
\liminf_{t\nearrow T^*}(T^*-t)\sqrt{|\ln(4(c\nu)^{-1}\|u(t)\|_{\dot H^{1/2}})|}\|u(t)\|_{\dot H^{5/2}}\geq c_1.\end{eqnarray}
\end{theorem}
The last result is the following
\begin{theorem}\label{theo3}
Let $u\in{\mathcal C}([0,T^*),H^{5/2}(\mathbb{R}^3)$ be a maximal solution of Navier-Stokes system. If, $T^*$ is finite then,  there is a universal constant $c_2>0$ such that
\begin{eqnarray}\label{enq113}
\liminf_{t\nearrow T^*}(T^*-t)\sqrt{|\ln(8\nu^{-1}\|u(t)\|_{\mathcal X^{-1})|}}\|u(t)\|_{\dot H^{5/2}}\geq c_2.
\end{eqnarray}
\end{theorem}

The paper is organized in the following way: In section $2$, we give some notations and important preliminary results. Section $3$ we prove the main result of this paper and we give some important remarks. The proof used standard Fourier techniques. In section $4$ and $5$ we give a proof respectively of theorem \ref{theo2} and \ref{theo3}. In section $6$, we give a simple proof of the explosion result in the space $\dot{H}^{3/2}$.\\

\section{\bf Notations and preliminary results}
\subsection{Notations}
In this section, we collect some notations and definitions that will be used later.\\
$\bullet$ The Fourier transformation is normalized as
$$
\mathcal{F}(f)(\xi)=\widehat{f}(\xi)=\int_{\mathbb R^3}\exp(-ix.\xi)f(x)dx,\,\,\,\xi=(\xi_1,\xi_2,\xi_3)\in\mathbb R^3.
$$
$\bullet$ The inverse Fourier formula is
$$
\mathcal{F}^{-1}(g)(x)=(2\pi)^{-3}\int_{\mathbb R^3}\exp(i\xi.x)g(\xi)d\xi,\,\,\,x=(x_1,x_2,x_3)\in\mathbb R^3.
$$
$\bullet$ The convolution product of a suitable pair of function $f$ and $g$ on $\mathbb R^3$ is given by
$$
(f\ast g)(x):=\int_{\mathbb R^3}f(y)g(x-y)dy.
$$
$\bullet$ If $f=(f_1,f_2,f_3)$ and $g=(g_1,g_2,g_3)$ are two vector fields, we set
$$
f\otimes g:=(g_1f,g_2f,g_3f),
$$
and
$$
{\rm div}\,(f\otimes g):=({\rm div}\,(g_1f),{\rm div}\,(g_2f),{\rm div}\,(g_3f)).
$$
$\bullet$ Let $(B,||.||)$, be a Banach space, $1\leq p \leq\infty$ and  $T>0$. We define $L^p_T(B)$ the space of all
measurable functions $[0,t]\ni t\mapsto f(t) \in B$ such that $t\mapsto||f(t)||\in L^p([0,T])$.\\
$\bullet$ The Sobolev space $H^s(\R^3)=\{f\in \mathcal S'(\R^3);\;(1+|\xi|^2)^{s/2}\widehat{f}\in L^2(\R^3)\}$.\\
$\bullet$ The homogeneous Sobolev space $\dot H^s(\R^3)=\{f\in \mathcal S'(\R^3);\;\widehat{f}\in L^1_{loc}\;and\;|\xi|^s\widehat{f}\in L^2(\R^3)\}$.\\
$\bullet$ The Lei-Lin space $\mathcal X^\sigma(\R^3)=\{f\in \mathcal S'(\R^3);\;\widehat{f}\in L^1_{loc}\;and\;|\xi|^\sigma\widehat{f}\in L^1(\R^3)\}$.\\
\subsection{Preliminary results}
In this section, we recall some classical results and we give new technical lemmas.\\
\begin{lem}\label{lem1}
We have $\mathcal X^{-1}(\R^3)\cap \mathcal X^{1}(\R^3)\hookrightarrow \mathcal X^0(\R^3)$. Precisely, we have
\begin{eqnarray}\label{enq1}
\|f\|_{\mathcal{X}^0(\mathbb{R}^3)}\leq\|f\|_{\mathcal{X}^{-1}(\mathbb{R}^3)}^{1/2}\|f\|_{\mathcal{X}^{1}(\mathbb{R}^3)}^{1/2},\;\forall f\in\mathcal X^{-1}(\R^3)\cap \mathcal X^{1}(\R^3).
\end{eqnarray}
\end{lem}
\noindent{\bf Proof.}
We can write
\begin{align*}
\|f\|_{\chi^0}&=\int |\widehat{f}|d\xi\\
              &\leq \int |\xi|^{1/2}|\widehat{f}|^{1/2}\frac{|\widehat{f}|^{1/2}}{|\xi|^{1/2}}d\xi.
\end{align*}
Cauchy Schwartz inequality gives the result.
\begin{lem}\label{lem2}
We have $\mathcal X^{-1}(\R^3)\cap \dot H^{5/2}(\R^3)\hookrightarrow \mathcal X^0(\R^3)$. Precisely, there is a constant $C_1>0$ such that
\begin{eqnarray}\label{enq2}
\|f\|_{\mathcal{X}^0(\mathbb{R}^3)}\leq C_1\|f\|_{\mathcal{X}^{-1}(\mathbb{R}^3)}^{1/2}\|f\|_{\dot H^{5/2}(\R^3)}^{1/2},\;\forall f\in\mathcal X^{-1}(\R^3)\cap \dot H^{5/2}(\R^3).
\end{eqnarray}
\end{lem}
\noindent{\bf Proof.}
For $R>0$ we have
$$\|f\|_{\chi^0}\leq \|f 1_{|D|<R}\|_{\chi^0}+\|f 1_{|D|>R}\|_{\chi^0}.$$
Cauchy Schwartz inequality gives
\begin{align*}
\|f 1_{|D|<R}\|_{\chi^0}&= \int_{|\xi|<R} |\widehat{f}(\xi)|d\xi\\
                        &\leq \int_{|\xi|<R} |\xi|\frac{|\widehat{f}(\xi)|}{|\xi|}d\xi\\
                        &\leq R \|f\|_{\chi^{-1}},
\end{align*}
and
\begin{align*}
\|f 1_{|D|>R}\|_{\chi^0}&= \int_{|\xi|>R} |\widehat{f}(\xi)|d\xi\\
                        &\leq \int_{|\xi|>R} |\xi|^{\frac{5}{2}}|\widehat{f}(\xi)||\xi|^{\frac{-5}{2}}d\xi\\
                        &\leq \big(\int_{|\xi|>R}|\xi|^{-5}\big)^{1/2}\|f\|_{\dot{H}^{5/2}}\\
                        &\leq \sqrt{\pi} R^{-1}\|f\|_{\dot{H}^{5/2}},
\end{align*}
To conclude, it suffices to take $R=\big(\frac{\|f\|_{\dot{H}^{5/2}}}{\|f\|_{\chi^{-1}}}\big)^{1/2}$

\begin{lem}\label{lem3}
We have $\dot H^{1/2}(\R^3)\cap  \mathcal X^{1}(\R^3)\hookrightarrow \mathcal X^0(\R^3)$. Precisely, there is a constant $C_2>0$ such that
\begin{eqnarray}\label{enq3}
\|f\|_{\mathcal{X}^0(\mathbb{R}^3)}\leq C_2\|f\|_{\dot H^{1/2}(\R^3)}^{1/2}\|f\|_{\mathcal{X}^{1}(\mathbb{R}^3)}^{1/2},\;\forall f\in\dot H^{1/2}(\R^3)\cap \mathcal X^{1}(\R^3).
\end{eqnarray}
\end{lem}
\noindent{\bf Proof.}
For $R>0$ we have
$$\|f\|_{\chi^0}\leq \|f 1_{|D|<R}\|_{\chi^0}+\|f 1_{|D|>R}\|_{\chi^0}.$$
Cauchy Schwartz inequality gives
\begin{align*}
\|f 1_{|D|<R}\|_{\chi^0}&= \int_{|\xi|<R} |\widehat{f}(\xi)|d\xi\\
                        &\leq \int_{|\xi|<R} |\xi|^{-1/2} |\xi|^{1/2}|\widehat{f}(\xi)|d\xi\\
                        &\leq \big(\int_{|\xi|<R} |\xi|^{-1}d\xi\big)^{1/2}\|f\|_{\dot{H}^{1/2}}\\
                        &\leq \sqrt{2\pi} R \|f\|_{\dot{H}^{1/2}}
\end{align*}
and
\begin{align*}
\|f 1_{|D|>R}\|_{\chi^0}&= \int_{|\xi|>R} |\widehat{f}(\xi)|d\xi\\
                        &\leq \int_{|\xi|>R} |\xi|^{-1}|\xi||\widehat{f}(\xi)|d\xi\\
                        &\leq R^{-1}\|f\|_{\chi^{1}}.
\end{align*}
To conclude, it suffices to take $R=\big(\frac{\|f\|_{\chi^1}}{\|f\|_{\dot{H}^{1/2}}}\big)^{1/2}$

\begin{lem}\label{lem4}
For $s\geq 0.$ If $f\in \dot{H}^s(\mathbb{R}^3)\cap \chi^1(\mathbb{R}^3)$, then, there exist a constant $C=C(s)>0$ such that
\begin{eqnarray}\label{enq4}
\big|\langle f.\nabla f, f\rangle_{\dot{H}^s}\big|\leq C \|f\|_{\chi^1}\|f\|_{\dot{H}^s}.
\end{eqnarray}
\end{lem}
\noindent{\bf Proof.}
If ${\rm div\,} f=0$ then, $\langle f.\nabla (|D|^s f), |D|^sf\rangle_{L^2}=0,$
which yields:
\begin{align*}
\big|\langle f.\nabla f, f\rangle_{\dot{H}^s}\big|&\leq \big| \langle |D|^{s}(f.\nabla f), |D|^{s} f\rangle_{L^2}\big|\\
                                                  &\leq \big| \langle |D|^{s}(f.\nabla f)-f.\nabla(|D|^{s}f), |D|^{s} f\rangle_{L^2}\big|\\
                                                  &\leq \||D|^{s}(f.\nabla f)-f.\nabla(|D|^{s}f)\|_{L^2} \||D|^{s} f\|_{L^2}\\
                                                  &\leq \||D|^{s}(f.\nabla f)-f.\nabla(|D|^{s}f)\|_{L^2} \|f\|_{\dot{H}^{s}}.
\end{align*}
we will estimate the first norm, we obtain:
\begin{align*}
\||D|^{s}(f.\nabla f)-f.\nabla(|D|^{s}f)\|_{L^2}&\leq \Big(\int \big||\xi|^{s}\mathcal{F}(f.\nabla f)-\mathcal{F}(f.\nabla(|D|^{s}f))\big|^2d\xi \Big)^{1/2}\\
                                                   &\leq \Big(\int \Big|\int|\xi|^{s}\widehat{f}(\xi-\eta)\widehat{\nabla f}(\eta)-\widehat{f}(\xi-\eta)|\eta|^{s}\widehat{\nabla f}(\eta)d\eta\Big|^2d\xi\Big)^{1/2}\\
                                                   &\leq \Big(\int \Big|\int\big(|\xi|^{s}-|\eta|^{s}\big)\widehat{f}(\xi-\eta)\widehat{\nabla f}(\eta)d\eta\Big|^2d\xi\Big)^{1/2}.
\end{align*}
 We have:
$$\Big||\xi|^{s}-|\eta|^s\Big|\leq s c^{s-1}|\xi-\eta|  \leq s|\xi-\eta| \big| |\xi|^{s-1}+|\eta|^{s-1} \big|,$$
$$|\xi|^{s-1}\leq 2^{s-1}|\xi-\eta|^{s-1}+2^{s-1}|\eta|^{s-1}\leq C_s( |\xi-\nu|^{s-1}+|\eta|^{s-1}).$$
Then, we get:
\begin{align*}
\||D|^{s}(f.\nabla f)-f.\nabla(|D|^{s})\|_{L^2}&\leq C\Big(\int \Big(\int|\xi-\eta||\xi-\eta|^{s-1}|\widehat{f}(\xi-\eta)|.|\widehat{\nabla f}(\eta)|d\eta\\
                                                &+\int|\xi-\eta||\widehat{f}(\xi-\eta)||\eta|^{s-1}|\widehat{\nabla f}(\eta)|d\eta\Big)^2d\xi\Big)^{1/2}\\
                                               &\leq C\Big(\int \Big(\int|\xi-\eta|^{s}|\widehat{f}(\xi-\eta)|\widehat{\nabla f}(\eta)|d\eta+\int|\widehat{\nabla f}(\xi-\eta)||\eta|^{s}|\widehat{ f}(\eta)|d\eta\Big|^2d\xi\Big)^{1/2}\\
                                               &\leq C\Big(\int \Big(\int|\xi-\eta|^{s}|\widehat{f}(\xi-\eta)||\widehat{\nabla f}(\eta)|d\eta\Big)^2d\xi\Big)^{1/2}\\
                                               &\leq C\|(|.|^s|\mathcal{F}(f)|)*|\mathcal{F}(\nabla f)|\|_{L^2},
\end{align*}
Young lemma yields,
\begin{align*}
\||D|^{s}(f.\nabla f)-f.\nabla(|D|^{s}f)\|_{L^2} &\leq C \|(|.|^s|\mathcal{F}(f)|)\|_{L^2}(\mathcal{F}(\nabla f))\|_{L^1}\\
                                                &\leq C \|f\|_{\dot{H}^s}\|f\|_{\chi^1}.
\end{align*}
Then, the proof is finished.

\begin{lem}\label{lem6}
For $f\in H^{7/2}(\R^3)$. There is a constant $C>0$ such that for all $0<\alpha<\beta<\infty$, we have
\begin{eqnarray}\label{enq61}
\|f\|_{\chi^1}\leq C\sqrt{4\pi} \Big[\alpha^{5/2}\|f\|_{L^2}+\sqrt{\ln(\frac{\beta}{\alpha})}\|f\|_{\dot{H}^{5/2}}+\frac{1}{\beta}\|f\|_{\dot{H}^{7/2}}\Big].
\end{eqnarray}
\begin{eqnarray}\label{enq62}
\|f\|_{\chi^1}\leq C\sqrt{4\pi} \Big[\alpha^2\|f\|_{\dot H^{1/2}}+\sqrt{\ln(\frac{\beta}{\alpha})}\|f\|_{\dot{H}^{5/2}}+\frac{1}{\beta}\|f\|_{\dot{H}^{7/2}}\Big].
\end{eqnarray}
\begin{eqnarray}\label{enq63}
\|f\|_{\chi^1}\leq C\sqrt{4\pi} \Big[\alpha^2\|f\|_{\mathcal X^{-1}}+\sqrt{\ln(\frac{\beta}{\alpha})}\|f\|_{\dot{H}^{5/2}}+\frac{1}{\beta}\|f\|_{\dot{H}^{7/2}}\Big].
\end{eqnarray}
\end{lem}
\noindent{\bf Proof.}

Let  $0<\alpha<\beta<\infty$ as:
$$
 \|f\|_{\chi^1}=\int |\xi||\widehat{f}(\xi)|d\xi= I_\alpha+J_{\alpha, \beta}+K_\beta,
$$
where: $$I_\alpha =\int_{|\xi|<\alpha}|\xi||\widehat{f}(\xi)|d\xi,$$
$$J_{\alpha,\beta} =\int_{\alpha<|\xi|<\beta}|\xi||\widehat{f}(\xi)|d\xi,$$
$$K_\beta =\int_{|\xi|>\beta}|\xi||\widehat{f}(\xi)|d\xi.$$
 Cauchy Schwartz inequality gives:
\begin{align*}
I_\alpha &=\int_{|\xi|<\alpha}|\xi||\widehat{f}(\xi)|d\xi\\
         &\leq \big(\int_{|\xi|<\alpha}|\xi|^2d\xi\big)^{1/2}\|f\|_{L^2}\\
         &\leq \sqrt{\frac{4\pi}{5}}\alpha^{5/2}\|f\|_{L^2},
\end{align*}
\begin{align*}
J_{\alpha,\beta} &=\int_{\alpha<|\xi|<\beta}|\xi||\widehat{f}(\xi)|d\xi\\
         &=\int_{\alpha<|\xi|<\beta}|\xi|^{-3/2}|\xi|^{5/2}|\widehat{f}(\xi)|d\xi\\
         &\leq \sqrt{4\pi} \big(\ln(\frac{\beta}{\alpha})\big)^{1/2}\|f\|_{\dot{H}^{5/2}},
\end{align*}
and,
\begin{align*}
K_\beta &=\int_{|\xi|>\beta}|\xi||\widehat{f}(\xi)|d\xi\\
         &=\int_{|\xi|>\beta}|\xi|^{-5/2}|\xi|^{7/2}|\widehat{f}(\xi)|d\xi\\
         &\leq \sqrt{4\pi} \Big(\frac{1}{2\beta^2}\Big)^{1/2}\|f\|_{\dot{H}^{7/2}}\\
         &\leq \sqrt{2\pi} \frac{1}{\beta}\|f\|_{\dot{H}^{7/2}}.
\end{align*}
Then, we can deduce (\ref{enq61}).
For the second estimate, we can write
\begin{align*}
I_\alpha &=\int_{|\xi|<\alpha}|\xi||\widehat{f}(\xi)|d\xi\\
         &=\int_{|\xi|<\alpha}|\xi|^{1/2}|\xi|^{1/2}|\widehat{u}(\xi)|d\xi\\
         &\leq (\int_{|\xi|<\alpha}|\xi|d\xi)^{1/2}\|f\|_{\dot{H}^{1/2}}\\
         &\leq \sqrt{4\pi} \big(\frac{\alpha^4}{4}\big)^{1/2}\|f\|_{\dot{H}^{1/2}}\\
         &\leq \sqrt{\pi} \alpha^2\|f\|_{\dot{H}^{1/2}},
\end{align*}
which gives the inequality (\ref{enq62}).
(\ref{enq63}) deduce from
 \begin{align*}
I_\alpha &=\int_{|\xi|<\alpha}|\xi||\widehat{f}(\xi)|d\xi\\
         &=\int_{|\xi|<\alpha}|\xi|^2|\xi|^{-1}|\widehat{f}(\xi)|d\xi\\
         &\leq  \alpha^2\|f\|_{\chi^{-1}},
\end{align*}
Then, the desired result is proved, and the proof of Lemma \ref{lem6} is finished.\\
\subsection{Remarks}
\begin{enumerate}
 \item[(i)]Leray showed that if the maximal data $T^*$ is finite then $$\|u(t)\|_{\dot{H}^1}\geq C(T^*-t)^{-1/4}$$
Interpolation inequality gives:
$$\|u(t)\|_{\dot{H}^1}\leq \|u^0\|^{3/5}_{L^2}\|u(t)\|^{2/5}_{\dot{H}^{5/2}}$$
then,
$$C(T^*-t)^{-1/4}\leq \|u^0\|^{3/5}_{L^2}\|u(t)\|^{2/5}_{\dot{H}^{5/2}}$$
which implies that
$$\lim_{t\rightarrow T^*} \|u(t)\|_{\dot{H}^{5/2}}=+\infty$$
if $T^*<\infty$ we can suppose that \begin{eqnarray}\label{R1}\|u(t)\|_{\dot{H}^{5/2}}>1\quad \forall t\in[0,T^*).\end{eqnarray}
\item[(ii)] if there is a time $t_0\in [0,T^*)$ such that $\|u(t_0)\|_{\dot{H}^{1/2}}<c\nu$ then $T^*=+\infty.$\\
Particulary, if $T^*<+\infty,$ then \begin{eqnarray}\label{R2}\|u\|_{\dot{H}^{1/2}}\geq c\nu\quad t\in [0,T^*).\end{eqnarray}
\item[(iii)] if there is a time $t_0\in [0,T^*)$ such that $\|u(t_0)\|_{\chi^{-1}}<\nu$ then $T^*=+\infty.$\\
Particulary, if $T^*<+\infty,$ then
\begin{eqnarray}\label{R3}\|u\|_{\chi^{-1}}\geq \nu\quad t\in [0,T^*).\end{eqnarray}
\end{enumerate}
\section{\bf Proof of Theorem\ref{theo1}}
 Let $u$ be a maximal solution of Navier-Stokes system in the space $C\big([0,T^*), \dot{H}^{5/2}(\mathbb{R}^3)\big)\cap L^2\big([0,T^*), \dot{H}^{7/2}(\mathbb{R}^3)\big).$\\
Suppose that  $u\in \dot{H}^{5/2}(\mathbb{R}^3)\cap \chi^1(\mathbb{R}^3)$.
Taking the norm of $\dot{H}^{5/2}(\mathbb{R}^3)$ and using lemma \ref{lem4}, we get:
\begin{align*}\partial_t \|u\|^2_{\dot{H}^{5/2}} +2\nu \|u\|^2_{\dot{H}^{7/2}} &\leq \big|\langle u.\nabla u, u\rangle_{\dot{H}^{5/2}}\big|\\
                                                                                &\leq C\|u\|_{\chi^1}\|u\|^2_{\dot{H}^{5/2}}
                                                                                \end{align*}
Using inequality (\ref{enq61}) we obtain:
\begin{align*}
\partial_t \|u\|^2_{\dot{H}^{5/2}} +2\nu \|u\|^2_{\dot{H}^{7/2}}&\leq C\|u\|^2_{\dot{H}^{5/2}}\|u\|_{\chi^1}\\
                                                                &\leq C \|u\|^2_{\dot{H}^{5/2}}\Big[\alpha^{2/5}\|u^0\|_{L^2}+\sqrt{\ln(\frac{\beta}{\alpha})}\|u\|_{\dot{H}^{5/2}}+\frac{1}{\beta}\|u\|_{\dot{H}^{7/2}}\Big]\\
                                                                &\leq C \|u\|^2_{\dot{H}^{5/2}}\Big[\alpha^{2/5}\|u^0\|_{L^2}+\sqrt{\ln(\frac{\beta}{\alpha})}\|u\|_{\dot{H}^{5/2}}\Big]+\frac{4\pi C\|u\|^4_{\dot{H}^{5/2}}}{\nu\beta^2}+\frac{\nu}{2}\|u\|^2_{\dot{H}^{7/2}}.
\end{align*}
Then, we have:
$$\partial_t \|u\|^2_{\dot{H}^{5/2}}\leq C\|u\|^2_{\dot{H}^{5/2}}\Big[\alpha^{5/2}\|u^0\|_{L^2}+\sqrt{\ln(\frac{\beta}{\alpha})}\|u\|_{\dot{H}^{5/2}}+\frac{\|u\|^2_{\dot{H}^{5/2}}}{\nu\beta^2}\Big].$$
Put $$t_0=\inf\{t\in[0,T^*),\; \|u(t)\|_{\dot{H}^{5/2}}= 2\|u^0\|_{\dot{H}^{5/2}}\},$$
then, we get:
$$4\|u^0\|^2_{\dot{H}^{5/2}}\leq \|u^0\|^2_{\dot{H}^{5/2}}+ C T^*\|u^0\|^2_{\dot{H}^{5/2}}\Big[\alpha^{5/2} \|u^0\|_{L^2}+\sqrt{\ln(\frac{\beta}{\alpha})}\|u^0\|_{\dot{H}^{5/2}}+\frac{\|u^0\|^2_{\dot{H}^{5/2}}}{\nu\beta^2}\Big].$$
By taking $\|u^0\|_{\dot{H}^{5/2}}>1$ by remark \ref{R1}, we can choose: $$\alpha = \|u^0\|^{2/5}_{\dot{H}^{5/2}}<\beta=\sqrt{\|u^0\|_{\dot{H}^{5/2}}},$$
 we obtain:
 \begin{align*}
 1&\leq C T^*\Big[ \|u^0\|_{\dot{H}^{5/2}}\|u^0\|_{L^2}+\sqrt{\ln(\|u^0\|_{\dot{H}^{5/2}})}\|u^0\|_{\dot{H}^{5/2}}+\frac{\|u^0\|_{\dot{H}^{5/2}}}{\nu}\Big]\\
  &\leq
  C T^*\|u^0\|_{\dot{H}^{5/2}}\Big[ \|u^0\|_{L^2}+\sqrt{\ln(\|u^0\|_{\dot{H}^{5/2}})}+\frac{1}{\nu}\Big].
 \end{align*}
 Let $t_1\in[0,T^*)$ be a instant such that,
 $$\sqrt{\ln(\|u(t)\|_{\dot{H}^{5/2}})}\geq 2(\|u^0\|_{L^2}+\frac{1}{\nu}),\quad \forall t\in [t_1,T^*)$$
 which yields: $$1\leq  C(T^*-t_1)\|u(t_1)\|_{\dot{H}^{5/2}}\sqrt{\ln(\|u(t_1)\|_{\dot{H}^{5/2}})}.$$
 By changing $t_1$ with any $t\in[t_1,T^*)$, we get the following estimate:
 $$\frac{1}{T^*-t}\leq C \|u(t)\|_{\dot{H}^{5/2}}\sqrt{\ln(\|u(t)\|_{\dot{H}^{5/2}})},\quad \forall t\in [t_1,T^*). $$
 Then, there is a constant $C_1>0$ for all $0\leq t<T^*$ we have:
 $$\frac{1}{C_1(T^*-t)\sqrt{|\ln(T^*-t)|}}\leq \|u(t)\|_{\dot{H}^{5/2}}.$$
 In fact:
  Let \begin{align*} f:&[4,+\infty[\rightarrow[4\sqrt{\ln(4)}, +\infty[\\
                                   &x\mapsto f(x)=x\sqrt{\ln(x)}\end{align*} be a continues bijection.\\
  Then, we have: $$y=x\sqrt{\ln(x)}\Longrightarrow\ln(y)=\ln(x)+ \ln(\sqrt{\ln(x)})\Longrightarrow\ln(y)\underset{x\rightarrow +\infty}{\sim} \ln(x)$$
  which implies
  $$x \underset{y\rightarrow +\infty}{\sim} \frac{y}{\sqrt{\ln(y)}}.$$

\section{\bf Blow up criterion in $\dot H^{5/2}$ with respect to $\dot H^{1/2}$ norm}
In this section we give a proof of Theorem \ref{theo2}.\\
Inequality (\ref{enq62}) implies
   $$\|u\|_{\chi^1}\leq C \Big[\alpha^2\|u\|_{\dot{H}^{1/2}}+\sqrt{\ln(\frac{\beta}{\alpha})}\|u\|_{\dot{H}^{5/2}}+\frac{1}{\beta}\|u\|_{\dot{H}^{7/2}}\Big].$$
  Then, we get
$$
\partial_t \|u\|^2_{\dot{H}^{5/2}} +2\nu \|u\|^2_{\dot{H}^{7/2}}\leq C\|u\|^2_{\dot{H}^{5/2}}\Big[\alpha^2\|u\|_{\dot{H}^{1/2}}+\sqrt{\ln(\frac{\beta}{\alpha})}\|u\|_{\dot{H}^{5/2}}\Big]+\frac{ C\|u\|^4_{\dot{H}^{5/2}}}{\nu\beta^2}+\frac{\nu}{2}\|u\|^2_{\dot{H}^{7/2}}.
$$
Moreover, we have:
$$\partial_t \|u\|^2_{\dot{H}^{5/2}}\leq  C\|u\|^2_{\dot{H}^{5/2}}\Big[\alpha^2\|u\|_{\dot{H}^{1/2}}+\sqrt{\ln(\frac{\beta}{\alpha})}\|u\|_{\dot{H}^{5/2}}+\frac{\|u\|^2_{\dot{H}^{5/2}}}{\nu\beta^2}\Big].$$
Let $t_0=\inf\{t\in[0,T^*),\; \|u(t)\|_{\dot{H}^{5/2}}= 2\|u^0\|_{\dot{H}^{5/2}}\},$ then we get:
$$\|u^0\|^2_{\dot{H}^{5/2}}\leq C T^*\|u^0\|^2_{\dot{H}^{5/2}}\Big[\alpha^2\sup_{0\leq s\leq t_0}\|u(s)\|_{\dot{H}^{1/2}}+\sqrt{\ln(\frac{\beta}{\alpha})}\|u^0\|_{\dot{H}^{5/2}}+\frac{\|u^0\|^2_{\dot{H}^{5/2}}}{\nu\beta^2}\Big].$$
Using remark (\ref{R2}), we can choose $\alpha$ and $\beta$ such that: $$\alpha = \sqrt{\frac{c\nu\|u^0\|_{\dot{H}^{5/2}}}{\sup_{0\leq t\leq s}\|u(t)\|_{\dot{H}^{1/2}}}}<\beta=\sqrt{\|u^0\|_{\dot{H}^{5/2}}}.$$
Then, we obtain:
$$
\|u^0\|^2_{\dot{H}^{5/2}}\leq C T^*\|u^0\|^2_{\dot{H}^{5/2}}\Big[\|u^0\|_{\dot{H}^{5/2}}+\sqrt{\ln\big(\frac{\sup_{0\leq t\leq s}\|u(t)\|_{\dot{H}^{1/2}}}{c\nu}\big)}\|u^0\|_{\dot{H}^{5/2}}+\frac{\|u^0\|_{\dot{H}^{5/2}}}{\nu}\Big].$$
Consequently, we get:
\begin{align*}
1&\leq  C T^*\|u^0\|_{\dot{H}^{5/2}}\Big[1+\sqrt{\ln\big(\frac{\sup_{0\leq t\leq s}\|u(t)\|_{\dot{H}^{1/2}}}{c\nu}\big)}+\frac{1}{\nu}\Big]\\
&\leq  C T^*\|u^0\|_{\dot{H}^{5/2}}\Big[C_\nu+\sqrt{\ln\big(\frac{\sup_{0\leq t\leq s}\|u(t)\|_{\dot{H}^{1/2}}}{c\nu}\big)}\Big].
\end{align*}
On the other hand, by interpolation, we have:
\begin{align*}\partial_t\|u\|^2_{\dot{H}^{1/2}}+2\nu \|u\|_{\dot{H}^{3/2}}&\leq \|u\otimes u\|_{\dot{H}^{1/2}}\|u\|_{\dot{H}^{3/2}}\\
                                                                          &\leq C_0\|u\|_{\dot{H}^{1/2}}\|u\|^2_{\dot{H}^{3/2}}\\
                                                                          &\leq C_0\|u\|^2_{\dot{H}^{1/2}}\|u\|_{\dot{H}^{5/2}}.
                                                                          \end{align*}
Then,  $$\frac{\frac{\partial_t\|u\|^2_{\dot{H}^{1/2}}}{c\nu}}{\frac{\|u\|^2_{\dot{H}^{1/2}}}{c\nu}}\leq C_0\|u\|_{\dot{H}^{5/2}}.$$
Integrating on $[0,t)$, we obtain:
$$\ln(\frac{\|u\|^2_{\dot{H}^{1/2}}}{c\nu})\leq \ln(\frac{\|u^0\|^2_{\dot{H}^{1/2}}}{c\nu})+C_0\int^t_0\|u(s)\|_{\dot{H}^{5/2}}ds.$$
Taking the sup over $[0,t_0]$, we get:
\begin{align*}
 \ln(\frac{\sup_{0\leq s\leq t_0}\|u(s)\|_{\dot{H}^{1/2}}}{c\nu})&\leq \ln(\frac{\|u^0\|^2_{\dot{H}^{1/2}}}{c\nu})+2C_0t_0\|u^0\|_{\dot{H}^{5/2}}\\
                                                    &\leq \ln(\frac{\|u^0\|^2_{\dot{H}^{1/2}}}{c\nu})+2C_0T^*\|u^0\|_{\dot{H}^{5/2}}.
\end{align*}
Therefore,
\begin{align*}
1&\leq  C_1 T^*\|u^0\|_{\dot{H}^{5/2}}\Big[C_\nu+\sqrt{\ln(\frac{\|u^0\|^2_{\dot{H}^{1/2}}}{c\nu})}+\sqrt{2T^*\|u^0\|_{\dot{H}^{5/2}}}\Big]\\
 &\leq  C_\nu T^*\|u^0\|_{\dot{H}^{5/2}}\Big[1+\sqrt{\ln(\frac{\|u^0\|^2_{\dot{H}^{1/2}}}{c\nu})}+2T^*\|u^0\|_{\dot{H}^{5/2}}\Big].
\end{align*}
Put $X=T^*\|u^0\|_{\dot{H}^{5/2}}$ and $a=1+\sqrt{\ln(\frac{\|u^0\|^2_{\dot{H}^{1/2}}}{c\nu})}$, we get:
$$X(a+X)\geq \frac{1}{ C_\nu}$$
$$P(X)=X^2+aX-\frac{1}{ C_\nu}\geq 0\;\Longrightarrow\; \Delta=a^2+\frac{1}{ C_\nu}>0 $$
The solution of $P$ are: $$\begin{cases}
                                        X_1=\frac{-a+\sqrt{\Delta}}{2}>0\\
                                        X_2=\frac{-a-\sqrt{\Delta}}{2}<0,
                          \end{cases}$$
then, we have:
$$X=T^*\|u^0\|_{\dot{H}^{5/2}}\geq X_1$$
which implies, $$ C_\nu T^*\|u^0\|_{\dot{H}^{5/2}}(1+a)\geq 1$$
$$  C_\nu T^*\|u^0\|_{\dot{H}^{5/2}}\sqrt{\ln(\frac{\|u^0\|^2_{\dot{H}^{1/2}}}{c\nu})}\geq 1. $$
We change the initial data with any $t\in [0,T^*)$, we obtain:
$$\liminf_{t\rightarrow T^*} (T^*-t)\|u(t)\|_{\dot{H}^{5/2}}\sqrt{\ln\big(\frac{\|u(t)\|^2_{\dot{H}^{1/2}}}{c\nu}\big)}\geq \frac{1}{ C_\nu }.$$

\section{\bf Blow up criterion in $\dot H^{5/2}$ with respect to $\mathcal X^{-1}$ norm}In this section we give a proof of Theorem \ref{theo3}.\\ Using (\ref{enq63}), we obtain:
$$\partial_t \|u\|^2_{\dot{H}^{5/2}}\leq C \|u\|^2_{\dot{H}^{5/2}}\Big[\alpha^2\|u\|_{\chi^{-1}}+\sqrt{\ln(\frac{\beta}{\alpha})}\|u\|_{\dot{H}^{5/2}}+\frac{\|u\|^2_{\dot{H}^{5/2}}}{\nu\beta^2}\Big].$$
Let  $$t_0=\inf\{t\in[0,T^*),\; \|u(t)\|_{\dot{H}^{5/2}}= 2\|u^0\|_{\dot{H}^{5/2}}\}.$$
which yields:
$$4\|u^0\|^2_{\dot{H}^{5/2}}\leq \|u^0\|^2_{\dot{H}^{5/2}}+ C T^*\|u^0\|^2_{\dot{H}^{5/2}}\Big[\alpha^2\sup_{0\leq s\leq t_0}\|u(s)\|_{\chi^{-1}}+\sqrt{\ln(\frac{\beta}{\alpha})}\|u^0\|_{\dot{H}^{5/2}}+\frac{\|u^0\|^2_{\dot{H}^{5/2}}}{\nu\beta^2}\Big]$$
Using remark (\ref{R3}) we can choose, $\alpha$ and $\beta$ such that: $$\alpha = \sqrt{\frac{c\nu\|u^0\|_{\dot{H}^{5/2}}}{\sup_{0\leq t\leq s}\|u(t)\|_{\chi^{-1}}}}<\beta=\sqrt{\|u^0\|_{\dot{H}^{5/2}}},$$
we get:
$$
\|u^0\|^2_{\dot{H}^{5/2}}\leq C T^*\|u^0\|^2_{\dot{H}^{5/2}}\Big[\|u^0\|_{\dot{H}^{5/2}}+\sqrt{\ln\big(\frac{\sup_{0\leq t\leq s}\|u(t)\|_{\chi^{-1}}}{c\nu}\big)}\|u^0\|_{\dot{H}^{5/2}}+\frac{\|u^0\|_{\dot{H}^{5/2}}}{\nu}\Big].$$
then, we have:
\begin{align*}
1&\leq C  T^*\|u^0\|_{\dot{H}^{5/2}}\Big[1+\sqrt{\ln\big(\frac{\sup_{0\leq t\leq s}\|u(t)\|_{\chi^{-1}}}{c\nu}\big)}+\frac{1}{\nu}\Big]\\
&\leq C_\nu  T^*\|u^0\|_{\dot{H}^{5/2}}\Big[1+\sqrt{\ln\big(\frac{\sup_{0\leq t\leq s}\|u(t)\|_{\chi^{-1}}}{c\nu}\big)}\Big].
\end{align*}
On the other hand, we have:
$$\|u(t)\|_{\chi^{-1}}\leq \|u^0\|_{\chi^{-1}}+\int^t_0\|u(s)\|_{\chi^{-1}}\|u(s)\|_{\dot{H}^{5/2}}ds.$$
By Granwall lemma, we obtain, for $0\leq t< T^*$:
 $$\|u(t)\|_{\chi^{-1}}\leq \|u^0\|_{\chi^{-1}}e^{\int^t_0\|u(s)\|_{\dot{H}^{5/2}}ds}.$$
By applying $\ln$ function, we get:
$$ \ln(\frac{\sup_{0\leq s\leq t_0}\|u(s)\|_{\chi^{-1}}}{c\nu})\leq \ln(\frac{\|u^0\|^2_{\chi^{-1}}}{c\nu})+2C't_0\|u^0\|_{\dot{H}^{5/2}},$$
and
\begin{equation}\label{E1}
 \ln(\frac{\sup_{0\leq s\leq t_0}\|u(s)\|_{\chi^{-1}}}{c\nu})\leq \ln(\frac{\|u^0\|^2_{\chi^{-1}}}{c\nu})+2C'T^*\|u^0\|_{\dot{H}^{5/2}}.
\end{equation}
by using (\ref{E1}), we can deduce:
$$
1\leq C_\nu T^*\|u^0\|_{\dot{H}^{5/2}}\Big[2+\sqrt{\ln(\frac{\|u^0\|_{\chi^{-1}}}{c\nu})}+2C'T^*\|u^0\|_{\dot{H}^{5/2}}\Big].
$$
then we have:
$$  C_\nu T^*\|u^0\|_{\dot{H}^{5/2}}\sqrt{\ln(\frac{\|u^0\|_{\chi^{-1}}}{c\nu})}\geq 1. $$
By changing the initial data with any $t\in [0,T^*)$, we obtain:
$$ \liminf_{t\rightarrow T^*}(T^*-t)\|u(t)\|_{\dot{H}^{5/2}}\sqrt{\ln\big(\frac{\|u(t)\|_{\chi^{-1}}}{c\nu}\big)}\geq \frac{1}{ C_\nu}.$$

\section{\bf General remarks} In this section we give a simple proof of the explosion result in $\dot{H}^{3/2}(\mathbb{R}^3)$, we give a simple proof of the following theorem given in \cite{5}

                    \begin{prop}
                             Let $ u\in C\big([0,T^*), (H^{3/2}(\mathbb{R}^3))^3\big)$ be a maximale solution of Navier-Stokes system. If $T^*$ is finite then,\\
                                      $$\|u(t)\|_{\dot{H}^{3/2}}\geq C\nu^{1/2}(T^*-t)^{-1/2}.$$
                     \end{prop}
\noindent{\it Proof.}  Taking the inner product in $\dot{H}^{3/2},$, we obtain:
                                   $$\langle\partial_tu,u\rangle_{\dot{H}^{3/2}}-\nu\langle\Delta u,u\rangle_{\dot{H}^{3/2}}+\langle u.\nabla u,u\rangle_{\dot{H}^{3/2}}=-\underbrace{\langle \nabla p,u\rangle_{\dot{H}^{3/2}}}_{0}.$$
                                   $\bullet$ we start by prove: $\langle u.\nabla u,u\rangle_{\dot{H}^{3/2}}\leq C \|u\|^2_{\dot{H}^{3/2}}\|u\|_{\dot{H}^{5/2}}.$\\
                                   we have:
                                   \begin{align*}
                                          \langle u.\nabla u,u\rangle_{\dot{H}^{3/2}} &= \int_{\mathbb{R}^3} |\xi|^3 \mathcal{F}(u.\nabla u)(\xi) \mathcal{F}(u)(-\xi)d\xi\\
                                          &=\int_{\xi}\int_{\eta} |\xi|^{3/2}\widehat{u}(\xi-\eta)\widehat{\nabla u}(\eta)|\xi|^{3/2}\widehat{u}(-\xi)d\eta d\xi\\
                                   \end{align*}
                                  By using, $$\langle f.\nabla g,g\rangle_{L^2}=0 {\quad\rm si\quad} {\rm div\;}f=0,$$ we get,
                                     $$\langle u.\nabla u,u\rangle_{\dot{H}^{3/2}}=\int_{\xi}\int_{\eta} |\xi|^{3/2}\widehat{u}(\xi-\eta)\widehat{\nabla u}(\eta)|\xi|^{3/2}\widehat{u}(-\xi)
                                          -|\eta|^{3/2}\widehat{u}(\xi-\eta)\widehat{\nabla u}(\eta)|\xi|^{3/2}\widehat{u}(-\xi)d\eta d\xi.$$
                                    Cauchy-Schawrtz inequality gives:
                                          $$\Big|\langle u.\nabla u,u\rangle_{\dot{H}^{3/2}}\Big|\leq \Big(\int_{\xi }\Big(\int_{\eta}(|\xi|^{3/2}-|\eta|^{3/2})|\widehat{u}(\xi-\eta)||\widehat{\nabla u}(\eta)|d\eta\Big)^2d\xi\Big)^{1/2}\|u\|_{\dot{H}^{3/2}}.$$
                                  For $\xi,\eta\in \mathbb{R}^3,$ we have\\
                                  $$\Big||\xi|^{3/2}-|\eta|^{3/2}\Big|\leq \frac{3}{2}\max(|\xi|,|\eta|)^{1/2}|\xi-\eta|\leq \frac{3}{2}(|\xi|^{1/2}+|\eta|^{1/2})|\xi-\eta|.$$
                                  Or we have $$|\xi|\leq |\xi-\eta|-|\eta|\leq 2 \max(|\xi-\eta|,|\eta|),$$
                                  then,$$|\xi|^{1/2}\leq \sqrt{2}(|\xi-\eta|^{1/2}+|\eta|^{1/2}),$$
                                  which yields,
                                  \begin{align*}
                                  \Big||\xi|^{3/2}-|\eta|^{3/2}\Big|&\leq \frac{3}{2}|\xi-\eta|(\sqrt{2}|\xi-\eta|^{1/2}+(1+\sqrt{2})|\eta|^{1/2})\\
                                  &\leq \frac{3(1+\sqrt{2})}{2}(|\xi-\eta|^{3/2}+|\eta|^{1/2}|\xi-\eta|).
                                  \end{align*}
                                  By using this inequality, we get:
                                  \begin{align*}
                                        &\Big(\int_{\xi}\Big(\int_{\eta}\Big||\xi|^{3/2}-|\eta|^{3/2}\Big||\widehat{u}(\xi-\eta)||\widehat{\nabla u}(\eta)|d\eta\Big)^2d\xi\Big)^{1/2}\\
                                        &\leq C\Big[\Big(\int_{\xi}\Big(\int_{\eta}|\xi-\eta|^{3/2}|\widehat{u}(\xi-\eta)||\widehat{\nabla u}(\eta)|d\eta\Big)^2d\xi\Big)^{1/2}\\
                                        &+\Big(\int_{\xi}\Big(\int_{\eta}|\xi-\eta||\widehat{u}(\xi-\eta)||\eta|^{1/2}|\widehat{\nabla u}(\eta)|d\eta\Big)^2d\xi\Big)^{1/2}\Big]\\
                                        &\leq I_1+I_2,
                                  \end{align*}
                                 where
                                  \begin{align*}
                                  \begin{cases}
                                       I_1&= \Big(\displaystyle\int_{\xi}\Big(\int_{\eta}|\xi-\eta|^{3/2}|\widehat{u}(\xi-\eta)||\widehat{\nabla u}(\eta)|d\eta\Big)^2d\xi\Big)^{1/2}= \|f_1g_1\|_{\dot{H}^0};\\
                                       f_1&= \mathcal{F}^{-1}(|\xi|^{3/2}|\widehat{u}(\xi)|)\\
                                       g_1&= \mathcal{F}^{-1}(|\widehat{\nabla u}(\xi)|),\\
                                       \end{cases}
                                  \end{align*}
                                 and
                                  \begin{align*}
                                       \begin{cases}
                                       I_2&=\Big(\displaystyle\int_{\xi}\Big(\int_{\eta}|\xi-\eta||\widehat{u}(\xi-\eta)||\eta|^{1/2}|\widehat{\nabla u}(\eta)|d\eta\Big)^2d\xi\Big)^{1/2}=\|f_2g_2\|_{\dot{H}^0};\\
                                       f_2&=\mathcal{F}^{-1}(|\xi||\widehat{u}(\xi)|)\\
                                       g_2&= \mathcal{F}^{-1}(|\xi|^{1/2}|\widehat{\nabla u}(\xi)|).
                                       \end{cases}
                                  \end{align*}
                                  Applying the product lower of  homogeneous Sobolev spaces, we obtain
                                  \begin{align*}
                                         I_1&\leq C\|f_1g_1\|_{\dot{H}^0}\\
                                         &\leq C \|f_1\|_{\dot{H}^1}\|g_1\|_{\dot{H}^{1/2}}\\
                                         &\leq C\|u\|_{\dot{H}^{5/2}}\|u\|_{\dot{H}^{3/2}}\\
                                 \end{align*}
                                 and
                                 \begin{align*}
                                         I_2&\leq C\|f_2g_2\|_{\dot{H}^0}\\
                                         &\leq C\|f_2\|_{\dot{H}^{1/2}}\|g_2\|_{\dot{H}^{1}}\\
                                         &\leq C\|u\|_{\dot{H}^{3/2}}\|u\|_{\dot{H}^{5/2}}.
                                  \end{align*}
                                  then, we get:
                                  $$\Big|\langle u.\nabla u,u\rangle_{\dot{H}^{3/2}}\Big|\leq C\|u\|^2_{\dot{H}^{3/2}}\|u\|_{\dot{H}^{5/2}}.$$
                                  which implied,
                                  $$ \partial_t\|u\|^2_{\dot{H}^{3/2}}+2\nu\|u\|^2_{\dot{H}^{5/2}}\leq C\|u\|^2_{\dot{H}^{3/2}}\|u\|_{\dot{H}^{5/2}}. $$
                                  Inequality $xy\leq \frac{x^2}{2}+\frac{y^2}{2}$, gives:
                                  \begin{align*}
                                         &\partial_t\|u\|^2_{\dot{H}^{3/2}}+2\nu \|u\|^2_{\dot{H}^{5/2}}\leq C\nu^{-1}\|u\|^4_{\dot{H}^{3/2}}+\frac{\nu}{2} \|u\|^2_{\dot{H}^{5/2}}\\
                                         \Longrightarrow\; & \partial_t\|u\|^2_{\dot{H}^{3/2}}\leq C\nu^{-1}\|u\|^4_{\dot{H}^{3/2}}.
                                  \end{align*}
                                  integrating over $[t, T^*)\subset[0, T^*)$, we get:
                                  $$ \|u(t)\|^2_{\dot{H}^{3/2}}\geq C\nu (T^*-t)^{-1}.$$
                                  Then, the proof of theorem \ref{theo2} is finished.


\begin{thebibliography}{10}

\bibitem{HB} H. Bahouri, J.Y Chemin and R. Danchin, {\it Fourier Analysis and Nonlinear Partial Differential
Equations\/}, Springer Verlag, 343p, 2011.

\bibitem{MC1} M. Cannone, {\it Harmonic analysis tools for solving the incompressible Navier-Stokes equations\/}, Diterot Editeur, Paris, 1995.


 599-624.


\bibitem{TK3}  H.Fujita,T.Kato, {\it On the Navier-Stokes initial value problem \/}, I.Arch.Ration.Mech.Anal.16 (1964) 269-315.

\bibitem{JB1} J. Benameur, {\it On the blow-up criterion of 3D Navier–Stokes equations.\/},Journal of Mathematical Analysis and Applications. Volume 371, Issue 2, 15 November 2010, Pages 719-727, (2010).
\bibitem{JB2} J. Benameur, {\it On the blow-up criterion of the periodic incompressible fluids.\/},Mathematical Methods in the Applied Sciences. Volume 36, Issue 2
30 January 2013
Pages 143–153, (2013).
\bibitem{JB3} J. Benameur, {\it Long Time Decay to the Lei-Lin solution of 3D Navier Stokes equation.\/} J.Math.Anal.Appl.(2015).



\bibitem{5} DS. McCormick, EJ. Olson, JC. Robinson, Jose L. Rodrigo, Alejandro Vidal-López, and Yi Zhou {\it Lower bounds on blowing-up solutions of the 3D Navier-Stokes equations in $\dot H^{3/2}$, $\dot H^{5/2}$, and $\dot B^{5/2}_{2,1}$.\/}, SIAM Journal on Mathematical Analysis. 2016, Vol. 48, No. 3, pp. 2119-2132, (2016).



\bibitem{CZ} A. Cheskidov and K. Zaya, {\it Lower bounds of potential blow-up solutions of three-dimentional Navier-Stokes equations in $\dot H^{3/2}$.\/} Journal of Mathematical Physics 57, 023101 (2016).


\bibitem{6} C. Cortissoz, J. A. Montero, and C. E. Pinilla , {\it On lower bounds for possible blow-up
solutions to the periodic Navier-Stokes equation\/}, J. Math. Phys., 55, 033101. (2014)

 \bibitem{4} J. C. Robinson, W. Sadowski, and R. P. Silva, {\it Lower bounds on blow up solutions of the three-dimensional Navier-Stokes equations in homogeneous Sobolev spaces\/}, J. Math. Phys, 53 (2012),115618.


\bibitem{TK1} T.Kato, {\it Quasi-Linear Equations of Evolution, With Application to Partial Differential Equations \/}, Lecture Notes in Math, vol.448,Sringer-Verlag,1975,pp. 25-70.

\bibitem{TK2} T.Kato. {\it $L^p$ -solution of the Navier Stokes in $\mathbb R^m$. With applications to weak solutions\/}, Math.Z.
187 (4)(1984) 471-480.

\bibitem{HKD} H.Koch.D.Tataru. {\it Well-posedness for the Navier Stokes equations \/}, Adv.Math.157(1)(2001) 22-35.

\bibitem{ZL} Z.Lei.F.Lin, {\it Global mild solutions of Navier Stokes equations  \/}, Comm.Pure Appl.Math.LXIV (2011) 1297 1304.

\bibitem{JL01} J.Leray. {\it Essai sur lr movement d'un liquide visqueux emplissant l'espace\/}, Acta Math.63 (1933) 22-25.

\bibitem{JL02} J.Leray. {\it Sur le movement d'un liquide visqueux emplissant l'espace\/}, Acta Math.63 (1) (1934) 193-248.



Indian University Mathematics Journal,vol.56,no. 3.pp.1157-1188,2007.

\end{thebibliography}
\end{document}